\documentclass{l4dc2026}

\title[Online Switching Control as Learning with Experts]{Online Switching Control as Learning with Experts}
\usepackage{times}

\usepackage{tikz}
\usepackage{comment}
\usepackage{wrapfig}

\makeatletter
\let\original@jmlrtps\ps@jmlrtps
\def\ps@jmlrtps{%
  \original@jmlrtps
  \let\@oddhead\@empty
  \let\@evenhead\@empty
}
\makeatother

\newcommand{\wtilde}[1]{\widetilde{#1}}
\newcommand{\calA}{{\mathcal A}} 
\newcommand{\calP}{{\mathcal P}} 
\newcommand{\bbR}{{\mathbb R}}
\newcommand{\bbE}{{\mathbb E}}
\newcommand{\bbP}{{\mathbb P}} 
\newcommand{\regret}{{R}} 
\newtheorem{assumption}{\textbf{Assumption}}
\newcommand{\ito}[1]{{\color{black}{#1}}}
\newcommand{\red}[1]{{\color{black}{#1}}}

\newcommand{\nakamura}[1]{{\color{black}{#1}}}

\coltauthor{\Name{Masataka Nakamura} \Email{mitsudomoe58massa@g.ecc.u-tokyo.ac.jp}\\
 \Name{Kaito Ito} \Email{kaito@g.ecc.u-tokyo.ac.jp}\\
 \Name{Hideaki Ishii} \Email{hideaki\_ishii@ipc.i.u-tokyo.ac.jp}\\
 \addr Department of Information Physics and Computing, The University of Tokyo, Tokyo, Japan}

\begin{document}

\maketitle

\begin{abstract}%
 This paper addresses the problem of online switching control for discrete-time nonlinear systems where a controller is sequentially selected from a finite set of candidates.
The goal is to generate an adaptive sequence of controllers to minimize the cumulative cost under unknown disturbances and cost functions.
Specifically, we consider a full-information feedback setting where realized disturbance effects and cost functions are available after each selection.
Using counterfactual costs computed from historical data,
we reduce the online switching control problem to the experts problem with switching costs. 
Then, we propose an adaptive switching algorithm based on the Hedge algorithm with a switching suppression mechanism. Notably, the proposed algorithm achieves an expected regret of $O(\sqrt{T\log N})$, where $T$ is the time horizon and $N$ is the number of candidate controllers, 
while ensuring uniform boundedness of the state.
This matches the optimal regret order for the standard experts problem. 
\end{abstract}


\section{Introduction}

\ito{In this paper, we address an adaptive switching control problem for discrete-time nonlinear systems under uncertain disturbances and control costs.
Adaptive switching control, also known as supervisory control, handles systems in uncertain environments by switching among multiple candidate controllers in response to incoming data~\citep{Morse2002supervisory,Hespanha2003hysteresis}.
Analogous problems have been tackled in the field of reinforcement learning as policy selection problems~\citep{Mazumdar2017multi,Li2018optimal}.}

\ito{Much effort has been devoted to ensuring stability under adaptive switching control~\citep{Vu2010supervisory}.
However, \red{existing optimality analyses in this framework have largely focused on asymptotic performance}, 
and the performance criteria are often tailored to specific objectives, such as stabilization~\citep{stefanovic2008safe,al2009switching}.
Consequently, \red{non-asymptotic performance guarantees for adaptive switching control under general cost functions remain largely unexplored.}
}

Recently, online control has attracted attention,
as it frames control problems under uncertainties as online learning problems
and provides non-asymptotic performance guarantees via regret analysis~\citep{agarwal2019online,Chekan2024learn,Muthirayan2025meta,Ito2025online,Watanabe2026}.
Regret is defined as the difference between the cumulative cost incurred by an online controller (or a learner) and that of the best controller chosen in hindsight.
This perspective has also been applied to switching control with a finite set of candidate controllers~\citep{anand2022online,li2023online} under the name of online switching control.
They investigated this problem 
under a bandit feedback setting, 
where only the incurred cost values are revealed to the learner. 
Then, they formulated the problem as an adversarial multi-armed bandit problem 
and proposed a batching-based switching algorithm that achieves a regret of
\ito{$\widetilde{O}(T^{2/3}N^{1/3})$}, where $T$ is the time horizon and $N$ is the number of controllers.

\red{Another important feedback model in online learning and control is the full-information setting}~\citep{Simchowitz2020improper,Orabona2023modern}. 
In this setting, the learner observes the cost function at each time step, rather than merely receiving the incurred cost value.
\red{Such a feedback setting arises, for example, when initially unknown cost parameters are revealed after each time step.}
Applications include adversarial target tracking~\citep{Abbasi2014tracking}, where the objective depends on the adversary's strategy, and industrial control, where the weight parameters of the cost function depend on realized electricity prices.
Recently, \citet{Xu2026online} proposed tracking a reference constructed from candidate controllers' counterfactual trajectories. Their method achieves sublinear regret for linear systems with convex costs, allowing inputs beyond the outputs of the candidate controllers.

In contrast, this paper investigates online switching control for nonlinear systems and nonconvex costs under full-information feedback, where the learner must select and apply one of the candidate controllers without modifying its output.
Our main contribution is to design an online switching algorithm that achieves a regret of $O(\sqrt{T\log N})$ while guaranteeing uniform boundedness of the state for a known nonlinear system. The regret bound holds under a uniform local Lipschitz-type condition on the cost functions and bounded disturbances.
The key idea is to combine learning with expert advice and stability-preserving switching among candidate controllers. Specifically, counterfactual costs computed from revealed cost functions and reconstructed disturbance effects allow us to reduce the learning component to the experts problem. By combining the Hedge algorithm with a sampling scheme that minimizes switching probabilities and a minimum dwell-time constraint, we achieve sublinear regret while ensuring state boundedness.
Unlike the \red{batching}-based approaches for bandit feedback in \citet{anand2022online,li2023online}, our method achieves $O(\sqrt{T\log N})$ regret without batching, matching the optimal rate
for the standard experts problem~\citep{geulenregret}.

\textbf{Organization:}
Section~\ref{sec:formulation} provides the problem formulation. 
In Section~\ref{sec:algorithm}, we propose an online switching algorithm based on a reduction to the experts problem. Section~\ref{sec:regret} establishes the regret guarantee. 
Section~\ref{sec:experiment} illustrates the effectiveness of the proposed method through numerical examples, followed by concluding remarks in Section~\ref{sec:conclusion}.

\textbf{Notation:}
Let $\mathbb{R}$ denote the set of real numbers.
The indicator function $\mathbb{I}[\cdot]$ equals $1$ if the condition inside the brackets is true, and $0$ otherwise.
We use $\widetilde{O}(\cdot)$ for the big-O notation hiding logarithmic factors.
Let $ \bbP[S] $ denote the probability of an event $ S $, and let $\mathbb{E}$ denote the expectation operator.
The Euclidean norm is denoted by $ \| \cdot \| $.

\section{Problem Formulation}\label{sec:formulation}
In this paper, we consider the online switching control problem 
for a discrete-time nonlinear system 
\begin{equation}\label{eq:system}
    x_{t+1} = f(x_t, u_t, w_t),   \qquad  \ito{t \ge 0},
\end{equation}
where $x_t \in \mathbb{R}^n$ is the state, 
$u_t\in \mathbb{R}^m$ is the control input, and $w_t\in\mathbb{R}^r$ is an unknown, \ito{nonstochastic\footnote{\ito{Although we can obtain a sublinear regret bound even for stochastic disturbances, this paper focuses on the nonstochastic case due to space limitations.}}} disturbance.
We assume that the state $x_t$ is fully observable, \ito{and the transition map \ito{$f : \bbR^n \times \bbR^m \times \bbR^r \rightarrow \bbR^n $} is known to the learner.\footnote{Knowledge of \(f\) is not essential; it suffices that the learner knows mappings \(\psi\) and \(F\) satisfying Assumption~\ref{asu:dy}.
}}

At each time step, the learner (also called the supervisor) chooses a controller from a finite set of candidate controllers:
\begin{equation*}
	\left\{ \pi_i : \bbR^n \rightarrow \bbR^m \mid i \in \calP :=  \{ 1,\ldots,N\} \right\} , 
\end{equation*}
where $ \pi_i $ is a state-feedback controller.
For the control input generated by the selected controller, the learner suffers a cost of $ c_t(x_t,u_t) $, where $ c_t : \mathbb{R}^n \times \mathbb{R}^m \rightarrow \mathbb{R} $ is a time-varying cost function. 
In this work, we consider the online switching control in a full-information setting. Specifically, the cost functions $ \{c_t\} $ are initially unknown, but after observing the state $ x_t $ and executing the control input $u_t$ at time $t$, the cost function $ c_t(\cdot, \cdot) $ is revealed to the learner. Consequently, the learner can use only the past cost functions $ \{c_s\}_{s=0}^{t-1} $.

Let $I_t \in \calP$ be the index of the selected controller at time $ t $ and $ T \ge 1 $ be the time horizon.
Our objective is to design an adaptive switching algorithm $\calA$ that minimizes the cumulative cost $ \sum_{t=0}^{T-1} c_t (x_t,u_t) $ under unknown disturbances and costs, where $ u_t = \pi_{I_t} (x_t) $.
We allow the switching algorithm to be randomized for the adaptation.
As the performance metric of $\calA$, we adopt the expected regret:
\begin{align}
    \regret_T(\mathcal{A}) &:= \ito{\bbE\! \left[\sum_{t=0}^{T-1} c_t(x_t(\mathcal{A}), u_t(\mathcal{A}))\right]} \!  - \min_{i\in\mathcal{P}} \sum_{t=0}^{T-1} c_t(x_t^{(i)},u_t^{(i)}) ,
\end{align}
where $ x_t(\calA) $ and $ u_t (\calA) $ (resp. $ x_t^{(i)} $ and $ u_t^{(i)} $) are the state and control input under algorithm $ \calA $ (resp. the fixed controller $ \pi_i $).
The best controller in hindsight is defined with full knowledge of all disturbances and cost functions $ \{(w_t, c_t) \}_{t=0}^{T-1} $.
Then, online switching control aims to minimize the regret.
A sublinear regret $\regret_T(\calA) = o(T)$ 
implies that the expected time-averaged cost is asymptotically no larger than that of 
the optimal controller 
in hindsight.
\ito{In what follows, for notational simplicity, we omit the argument $ \calA $.}

\ito{In the rest of this section, we introduce assumptions for establishing a sublinear regret of the online switching control problem.}
\ito{In this paper, the disturbances and the cost functions may be chosen adversarially, subject to the following assumptions.}
\begin{assumption}[Oblivious adversary]\label{adverse}
    The sequences $\{w_t\}_{t=0}^{T-1}$ and $\{c_t\}_{t=0}^{T-1}$ are fixed before the first control input $u_0$ is selected, independently of the learner's randomization.
\end{assumption}

\begin{assumption}[Disturbances]\label{ass:dis}
    There exists a known constant $W\geq 0$ such that $\|w_t\|\leq W$ for all $t\ge 0$.
\end{assumption}

The following local Lipschitz continuity assumption is also used in \citep{li2023online,kim2024online} and holds, e.g., for quadratic tracking costs with bounded references and weights.

\begin{assumption}[Cost functions]\label{ass:cost}
There exist $L_{c1}, L_{c2}\ge 0$ such that for all
$t\ge0$, $x,x'\in\mathbb{R}^n$ and $u,u'\in\mathbb{R}^m$,
the cost function $c_t$ satisfies
\begin{equation*}
    |c_t(x,u)-c_t(x', u')|\leq \left[L_{c1}\{\max(\|x\|, \|x'\|)+\max(\|u\|, \|u'\|)\}+L_{c2} \right](\|x-x'\|+\|u-u'\|).
\end{equation*}
\end{assumption}

Following \citet{li2023online}, we assume exponential input-to-state stability (E-ISS) and incremental exponential stability ($\delta$-ES) for the candidate controllers.
\begin{definition}[E-ISS]
    Under a controller $\pi$, the system $x_{t+1} = f(x_t, \pi(x_t),w_t)$ is said to 
    be E-ISS with parameters $ \kappa \ge 1 $, $ \rho \in (0,1) $, and $ \beta \ge 0 $
    if, for any initial state $x_0\in \mathbb{R}^n$ and any disturbance sequence $\{w_t\}$ satisfying $\|w_t\|\leq W$ for all $t\geq 0$, the following inequality holds:
    \begin{equation}\label{eq:E-ISS}
        \|x_t\|\leq \kappa\rho^t\|x_0\| + \beta W,\qquad \forall t\ge 0.
    \end{equation}
\end{definition}

For an initial state $\xi$,
let $x_t(\xi)$ denote the state trajectory of the system
$x_{t+1}=f(x_t,\pi(x_t),w_t)$ starting from the initial
state $x_0=\xi$ under the controller $\pi$.
\begin{definition}[$\delta$-ES]
    \nakamura{Under a controller $\pi$,} the system $x_{t+1} = f(x_t, \pi(x_t), w_t)$ is said to be $\delta$-ES with parameters $ \kappa_\delta \ge 1 $ and $ \rho_\delta \in (0,1) $ if, for any $\xi^1, \xi^2\in \mathbb{R}^n$ and any disturbance sequence $\{w_t\}$ satisfying $\|w_t\|\leq W$ for all $t\geq 0$, the following inequality holds:
    \begin{equation}
            \left\| x_t(\xi^1) - x_t(\xi^2) \right\| \le \kappa_\delta \rho_{\delta}^t \! \left\|\xi^1-\xi^2 \right\|,\qquad \forall t\ge 0.
    \end{equation}
\end{definition}

\ito{By abuse of terminology}, we \ito{say that a controller $\pi$ is} E-ISS (resp. $\delta$-ES) when its
closed-loop system is E-ISS (resp. $ \delta $-ES).
While prior work on unknown systems allows potentially destabilizing candidate controllers~\citep{li2023online}, we consider known dynamics, where stabilizing controllers can be selected in advance, eliminating the need to include destabilizing candidates. 
\begin{assumption}[Candidate controllers]\label{ass:controller}
    All candidate controllers $\pi_1,\ldots,\pi_N$ satisfy $(\kappa, \rho, \beta)$-E-ISS and $(\kappa_\delta, \rho_\delta)$-$\delta$-ES with common known parameters.
    Moreover, each controller $\pi_i$ is Lipschitz continuous with common $L_{\pi}>0$, i.e., for all $x, x'\in\mathbb{R}^n$, we have $\|\pi_i(x)-\pi_i(x')\|\leq L_{\pi}\|x-x'\|$.
\end{assumption}

E-ISS ensures boundedness under each fixed controller, while $\delta$-ES is utilized to ensure that the cumulative cost penalty induced by each switch remains bounded through exponential forgetting of initial conditions.
\red{For linear time-invariant systems, Schur-stabilizing state-feedback controllers satisfy both properties.}
Lastly, let $\bar{\pi}_0 > 0$ satisfy $\max_{i\in\mathcal{P}}\|\pi_i(0)\|\leq \bar{\pi}_0$.


\red{Our approach is to reduce the online switching control problem to the experts problem, where the loss of each expert is revealed after every time step. In our setting, however, observing the cost function \(c_t\) alone is not sufficient to compute the losses of all candidate controllers, since their state trajectories depend on the unknown disturbances. Therefore, we impose the following assumption, which allows us to reproduce the effects of the realized disturbances without necessarily recovering the disturbances themselves.}
\begin{assumption}[Replayability of disturbance effects]\label{asu:dy}
    For the transition map $f$, there exist an auxiliary space $ \mathcal{Z} $ and known mappings
    $ \psi: \bbR^n \times \bbR^m \times \bbR^n \to \mathcal{Z} $ and
    $ F: \bbR^n \times \bbR^m \times \mathcal{Z} \to \bbR^n $ such that, for any $ x, \tilde{x} \in \bbR^n $, $ u, \tilde{u} \in \bbR^m $, and $ w \in \bbR^r $, 
    \begin{equation}\label{eq:replay}
        F \! \left( \tilde{x}, \tilde{u}, \psi(x, u, f(x,u,w)) \right)
        = f(\tilde{x}, \tilde{u}, w).
    \end{equation}
\end{assumption}

The mapping $ \psi $ extracts the effect of the realized disturbance, and $ F $ replays it. 
For additive disturbances $ f(x,u,w) = g(x,u) + Gw $, one can use $ \psi(x,u,x^+) := x^+ - g(x,u) $ and
$ F(\tilde{x},\tilde{u},z) := g(\tilde{x},\tilde{u}) + z $.
The assumption also covers a state-dependent
disturbance channel $ f(x,u,w) = g(x,u) + G(x,u)w $ with $ G(x,u) $ of full
column rank, for which $ \psi(x,u,x^+) := G(x,u)^{\dagger}(x^+ - g(x,u)) $
and $ F(\tilde{x},\tilde{u},z) := g(\tilde{x},\tilde{u}) + G(\tilde{x},
\tilde{u}) z $ verify \eqref{eq:replay}, where $ \dagger $ denotes the
Moore--Penrose pseudoinverse.

\section{Algorithm Design}\label{sec:algorithm}
In this section, we present our online switching algorithm.

\subsection{Reduction to the Experts Problem}
\red{Under Assumption~\ref{asu:dy}, the observed state transitions allow us to reconstruct the effect of the realized disturbances. Together with the revealed cost functions, this enables us to simulate the trajectories and costs associated with each fixed candidate controller.}
\begin{definition}[Counterfactual states, inputs, and costs] \label{def:counterfactual}
    For each controller $i \in \mathcal{P}$ and time $t \ge 0$, we define the counterfactual state, input, and cost $ (x_t^{(i)}, u_t^{(i)}, c_t^{(i)} ) \in \mathbb{R}^n \times \bbR^m \times \bbR $ as follows:
    \begin{align*}
        &x_0^{(i)}=x_0,\quad x_{t+1}^{(i)}=f(x_t^{(i)}, \pi_i(x_t^{(i)}),w_t),\quad u_t^{(i)} := \pi_i(x_t^{(i)}),\quad c_t^{(i)}:= c_t(x_t^{(i)},u_t^{(i)}).
    \end{align*}
\end{definition}

In our setting, these quantities can be computed.
Indeed, from the observed transition $ (x_t, u_t, x_{t+1}) $, the effect of
$ w_t $ is extracted as
$
    z_t := \psi(x_t, u_t, x_{t+1})
$,
and then \eqref{eq:replay} yields, for all $ i \in \calP $,
\begin{equation}\label{eq:counterfactual_update}
    F ( x_t^{(i)}, \pi_i(x_t^{(i)}), z_t )
    = f ( x_t^{(i)}, \pi_i(x_t^{(i)}), w_t )
    = x_{t+1}^{(i)}.
\end{equation}
Hence, the counterfactual states are computed recursively from
$ x_0^{(i)} = x_0 $.
Moreover, since the cost function $ c_t $ is revealed in the full-information
setting, the counterfactual costs $ c_t^{(i)} $ are also computed.
When determining the control input $ u_t $ at time $ t $, the counterfactual
costs $ \{c_s^{(1)},\ldots,c_s^{(N)} \}_{s=0}^{t-1} $ are thus available.

We regard each controller $i\in \calP$ as an expert with loss \red{$c_t^{(i)}$} and apply the Hedge algorithm~\citep{freund1997decision}.
Specifically, Hedge selects expert $i$ with probability
$p_{t,i}$ at time $t$ by exponential weighting of each expert's cumulative loss with learning rate $ \eta > 0 $:
\begin{equation}\label{hedge:update}
    p_{t,i}:=\frac{\exp\bigl(-\eta\sum_{s=0}^{t-1}c_s^{(i)}\bigr)}
    {\sum_{k=1}^N\exp\bigl(-\eta\sum_{s=0}^{t-1}c_s^{(k)}\bigr)},
    \qquad p_{0,i}=\frac1N.
\end{equation}
Under this scheme, experts with lower cumulative losses are 
selected with higher probabilities, 
while other experts retain a chance of being selected, avoiding premature commitment to a single expert.

By Assumption~\ref{adverse}, the counterfactual costs do not depend
on the learner's selections. Moreover, the counterfactual trajectories and
the range of the costs are uniformly bounded as follows.
The proof is given in Appendix~\ref{app:vbo}.
\begin{lemma}\label{lem:vbo}
Suppose that Assumptions~\ref{ass:dis}--\ref{ass:controller} hold.
Then, for any $ t \ge 0 $ and $ i \in \calP $, it holds that
\begin{align}
    \| x_t^{(i)} \| \le X_{\mathrm{c}}, \qquad \| u_t^{(i)} \| \le U_{\mathrm{c}},
    \label{eq:Xv_Uv}
\end{align}
where $ X_{\mathrm{c}} := \kappa \| x_0 \| + \beta W $ and
$ U_{\mathrm{c}} := L_\pi X_{\mathrm{c}} + \bar{\pi}_0 $.
Moreover, for any $ t \ge 0 $ and $ i, j \in \calP $, we have
\begin{align}
    |c_t^{(i)} - c_t^{(j)}| \le \Delta_{\mathrm{c}}, \label{eq:Rc}
\end{align}
where $ \Delta_{\mathrm{c}} := 2 \left\{ L_{c1}(X_{\mathrm{c}} + U_{\mathrm{c}}) + L_{c2} \right\}
(X_{\mathrm{c}} + U_{\mathrm{c}}) $.
\end{lemma}

Let $m_t:=\min_{i\in\calP}c_t^{(i)}$ and
$M_t:=\max_{i\in\calP}c_t^{(i)}$.
Then $M_t-m_t\le\Delta_{\mathrm c}$.
We use the following straightforward extension of the
standard Hedge regret bound for losses in $[0,1]$
\citep{cesa2006prediction} to general bounded loss ranges.
\begin{proposition}[Regret of Hedge]\label{pro:hedge}
Suppose that Assumptions~\ref{ass:dis}--\ref{ass:controller} hold.
Let $p_{t,i}$ be defined by \eqref{hedge:update} and $\Delta_{\mathrm c}$ be as in Lemma~\ref{lem:vbo}. Then, for any $\eta > 0$, it holds that
\begin{equation}\label{eq:hedge_range}
\sum_{t=0}^{T-1} \sum_{i\in \calP} p_{t,i} c_t^{(i)} - \min_{i\in \calP} \sum_{\ito{t=0}}^{T-1} c_t^{(i)}
\leq \frac{\log N}{\eta}+\frac{\eta T\Delta_{\mathrm c}^2}{8} .
\end{equation}
\end{proposition}

\ito{Typically, Hedge performs an independent selection of an expert at each time step 
according to the distribution $ \{p_{t,i} \}_{i=1}^N $.
However, this independence can cause frequent switching.
Such switches can destabilize the system, 
even if each controller is stabilizing~\citep{liberzon2003switching}.
Moreover, controller switching can cause a discrepancy
between the actual cost $c_t(x_t,u_t)$ and the
counterfactual cost $c_t^{(I_t)}$, resulting in
a switching penalty in the regret analysis,
as shown in the proof of Theorem~\ref{theorem: regret}.}

To address these issues, we propose a two-stage
controller selection mechanism.

\subsection{Lazy Sampling}

\red{As the first stage}, we employ lazy sampling
based on maximal coupling~\citep{Sherman2021lazy} (Algorithm~\ref{alg:Lazy sampling})
to generate a tentative selection $J_t\in\mathcal{P}$,
which is then passed to the dwell-time filter to determine \red{the selection $I_t$ to be applied.}
This method preserves the Hedge marginal $J_t\sim p_t$
while achieving the minimum switching probability
$\mathbb{P}[J_{t+1}\neq J_{t}] = \frac{1}{2}\sum_{i\in\mathcal P}|p_{t+1,i}-p_{t,i}|$.
The following proposition establishes these properties.
The proof is given in Appendix~\ref{app:theorem:switch}.

\begin{algorithm2e}[tb]
\footnotesize
\caption{Lazy sampling (maximal coupling)}
\label{alg:Lazy sampling}

\KwIn{Probability distributions $p_{t-1}, p_t$ on $\mathcal{P}$, and a selection $J_{t-1}\in\mathcal{P}$.}

Define $k := J_{t-1}$.

With probability
$\min\!\left\{1,p_{t,k}/p_{t-1,k}\right\}$,
maintain the selection: set $J_t \leftarrow k$.

Otherwise, sample $J_t\in\mathcal{P}$ from the distribution
$q_t$ defined by $q_{t,i}\propto \max\{0,\,p_{t,i}-p_{t-1,i}\}$.

\end{algorithm2e}

\begin{proposition}[Lazy sampling]\label{theorem:switch}
Suppose that Assumptions~\ref{adverse}--\ref{ass:controller} hold.
Let $J_0\sim p_0$, and generate $\{J_t\}_{t=1}^{T-1}$ by
Algorithm~\ref{alg:Lazy sampling} with Hedge and a constant learning rate $\eta>0$.
Then, for all $i\in\calP$ and $t\in\{0,\ldots,T-1\}$, we have
$\mathbb P[J_t=i]=p_{t,i}$, and
\begin{equation}\label{switchprobability}
    \mathbb P[J_{t+1}\ne J_t] = \frac{1}{2} \sum_{i\in \calP} | p_{t+1,i} - p_{t,i} | \le\frac{\eta \Delta_{\mathrm c}}4,
    \qquad t\in\{0,\ldots,T-2\}.
\end{equation}
Furthermore, the total number of switches
$N_J:=\sum_{t=1}^{T-1}\mathbb I[J_{t}\ne J_{t-1}]$ satisfies
    $\mathbb E[N_J]\le\frac{\eta \Delta_{\mathrm c}(T-1)}4$.
\end{proposition}

\subsection{Dwell-Time Filter}\label{dtf}
\red{Although lazy sampling minimizes the switching probability, it does not guarantee a minimum interval between consecutive switches. Consequently, the state may still become unbounded.} 
To guarantee the boundedness of the state, \red{as the second stage in the controller,} we introduce a deterministic filter, which we call the dwell-time filter.
This enforces a minimum dwell time $\tau$, 
requiring the system 
to use the same controller for at least $\tau$ time steps.
Initialize $I_0 = J_0 $ and $t_{\mathrm{last}} = 0$.
For $t \ge 1$, let $t_{\mathrm{last}}$ denote the last time when the controller was switched before the current time $ t $, or $0$ if no switch has occurred.
Then, $I_t$ is determined based on the elapsed time $t - t_{\mathrm{last}}$ as follows: 
\begin{align}
            \ito{I_t = 
    \begin{cases}
    I_{t-1} (=I_{t_{\mathrm{last}}}) &\text{if} \ t-t_{\mathrm{last}}<\tau \ \ (\text{locked period}), \\
    J_t &\text{if} \ t-t_{\mathrm{last}}\geq\tau \ \ (\text{unlocked period}) . \label{eq:dwell_time_filter}
    \end{cases}}
        \end{align}
If $I_t\neq I_{t-1}$, the controller switches to $ \pi_{I_t} $ at time $ t $, and we update $t_{\mathrm{last}}\gets t$.
By construction, the total number of switches $N_I:=\sum_{t=1}^{T-1} \mathbb{I}[I_t\neq I_{t-1}]$ generated by the dwell-time filter satisfies
$
    N_I\leq N_J
$.

By appropriately setting the minimum dwell time $\tau$ based on the parameters of the E-ISS $(\kappa, \rho, \beta)$, 
the uniform boundedness of the state $x_t$ is guaranteed. The proof is deferred to Appendix~\ref{app:statebound}.
  \begin{proposition}[Uniform boundedness of the state]\label{th:statebound}
    Suppose that Assumptions \ref{ass:dis} and \ref{ass:controller} hold.
    Choose a dwell time $\tau$ such that $\tau>\frac{\log \kappa}{\log (1/\rho)}$. 
    Then, for any tentative sequence $\{J_t\}_{t=0}^{T-1} \in \calP^{T} $, the controller sequence $\{\pi_{I_t}\}_{t=0}^{T-1}$ generated by the dwell-time filter~\eqref{eq:dwell_time_filter} satisfies
    \begin{equation}
        \|x_t\| \leq \kappa\|x_0\| + \left(\frac{\kappa}{1-\kappa\rho^\tau} + 1\right)\beta W =: X,\qquad \forall t \in \{0,\ldots,T\} .\label{eq:X_def} 
    \end{equation}
\end{proposition}

In summary, the proposed algorithm consists of  
three components: 
updating a probability distribution over the controllers using Hedge, 
reducing the switching probability via lazy sampling, 
and enforcing a minimum dwell time through the dwell-time filter. 
The overall procedure is summarized in Algorithm \ref{alg:proposed_algorithm}.

\begin{algorithm2e}[tb]
\footnotesize
\caption{Hedge with lazy sampling and dwell-time filter}
\label{alg:proposed_algorithm}

\KwIn{Horizon $T \ge 1$, dwell time $\tau\ge1$, learning rate $\eta>0$.}

Initialize $p_{0,i}\gets 1/N$ and
$x_0^{(i)}\gets x_0$ for all $i\in\mathcal P$.

\For{$t=0,\ldots,T-1$}{

    \eIf{$t=0$}{
        Sample $J_0\sim p_0$ and set $I_0\gets J_0$.
    }{
        Generate $J_t$ by lazy sampling
        (Algorithm~\ref{alg:Lazy sampling}).

        Determine $I_t$ by the dwell-time filter
        \eqref{eq:dwell_time_filter}.
    }

    Apply $u_t=\pi_{I_t}(x_t)$ to the system.

    Receive $c_t$ and incur cost $c_t(x_t,u_t)$.

    Compute
    $c_t^{(i)}\gets
    c_t(x_t^{(i)},\pi_i(x_t^{(i)}))$
    for all $i\in\mathcal P$.

    Update $p_{t+1,i}\propto
    p_{t,i}e^{-\eta c_t^{(i)}}$
    for all $i\in\mathcal P$, with $\sum_i p_{t+1,i}=1$.

    Observe $x_{t+1}$ and compute
    $z_t\gets\psi(x_t,u_t,x_{t+1})$.

    Update
    $x_{t+1}^{(i)}
    \gets F(x_t^{(i)},\pi_i(x_t^{(i)}),z_t)$
    for all $i\in\mathcal P$.
}
\end{algorithm2e}

\section{Regret Guarantee}\label{sec:regret}
In this section, as the main result, 
we provide an upper bound on the regret of the proposed online switching control (Algorithm \ref{alg:proposed_algorithm}).

\begin{theorem}\label{theorem: regret}
    Suppose that \ito{Assumptions \ref{adverse}--\ref{asu:dy}} hold. 
    Choose a constant learning rate $\eta>0$ and a dwell time $\tau$ such that $\tau > \frac{\log \kappa}{\log (1/\rho)}$.
    Then, the expected regret of Algorithm~\ref{alg:proposed_algorithm} satisfies
    \begin{align}
            \regret_T \leq \frac{\log N}{\eta}
    + \frac{\eta T \Delta_{\mathrm{c}}}{4} \left\{ \alpha + \Delta_{\mathrm{c}}\left(\tau + \frac{1}{2}\right)\right\} .
        \label{regretstandard}
    \end{align}
    Here, $ \Delta_{\mathrm{c}} $ is defined in Lemma~\ref{lem:vbo}, and $\alpha$ is defined as
    \begin{align}
        \alpha :=
        \frac{\left\{L_{c1}\big((1+L_{\pi})X + \bar{\pi}_0\big) + L_{c2}\right\} \!
        (1+L_{\pi})(X+X_{\mathrm{c}})\kappa_{\delta}}{1-\rho_{\delta}}, \label{eq:alpha}
    \end{align}
    \ito{where $ X_{\mathrm{c}} $ and $ X $ are defined in \eqref{eq:Xv_Uv} and \eqref{eq:X_def}, respectively.}
    Moreover, setting $\eta = \Theta\!\left(\sqrt{T^{-1}\log N}\right)$ to minimize the upper bound in \eqref{regretstandard} yields
    $
        \regret_T = O\!\left(\sqrt{T \log N}\right)
    $.
\end{theorem}

\begin{proof}
    The proof relies on decomposing the regret into the difference between the actual and counterfactual costs, the cost discrepancy introduced by the dwell-time filter, and the Hedge regret:
\begin{align}
\regret_T=\mathbb E\!\left[\sum_{t=0}^{T-1}
 \big(c_t(x_t,u_t)-c_t^{(I_t)}\big)\right]
+\mathbb E\!\left[\sum_{t=0}^{T-1}
 \big(c_t^{(I_t)}-c_t^{(J_t)}\big)\right]+\mathbb E\!\left[\sum_{t=0}^{T-1} c_t^{(J_t)}\right]
-\min_{i\in\calP}\sum_{t=0}^{T-1} c_t^{(i)}.
\label{decomposition}
\end{align}

First, we bound the first term. Let $0=s_0<s_1<\cdots<s_{N_I}<T$ be the actual switching times, where $s_0$ represents initialization, and set $s_{N_I+1}:=T$. For $s_k\le t<s_{k+1}$, the actual and counterfactual trajectories evolve under the same controller and disturbance sequence. Thus, $\delta$-ES and the state bounds in Lemma~\ref{lem:vbo} and Proposition~\ref{th:statebound} give
$
\|x_t-x_t^{(I_t)}\|
\le\kappa_\delta\rho_\delta^{t-s_k}
\|x_{s_k}-x_{s_k}^{(I_t)}\|
\le\kappa_\delta\rho_\delta^{t-s_k}(X+X_c).
$
Regarding the control input, it follows from the Lipschitz continuity of the controller that
$
    \|u_{t} - u_{t}^{(I_{t})}\| \leq  L_{\pi}\|x_{t} - x_{t}^{(I_{t})}\| 
$.
Therefore, by using the local Lipschitz continuity of the cost function, the difference between the actual and counterfactual costs is bounded as
\begin{align*}
    &c_t(x_t, u_t) -c_t(x_t^{(I_t)}, u_t^{(I_t)})\\
    &\leq \left[L_{c1} \! \left\{\max(\|x_t\|,\|x_t^{(I_t)}\|) + \max(\|u_t\|,\|u_t^{(I_t)}\|)\right\} + L_{c2}\right](\|x_t - x_t^{(I_t)}\| + \|u_t - u_t^{(I_t)}\|)\\
    &\leq \left\{L_{c1}(X + L_{\pi}X + \bar{\pi}_0) + L_{c2}\right\}
    (1+L_{\pi})\kappa_\delta \rho_{\delta}^{t-s_k}(X+X_{\mathrm{c}}) =(1-\rho_{\delta})\alpha \rho_{\delta}^{t-s_k}.
\end{align*}
Summing this difference within the interval $ \{s_k,\ldots,s_{k+1}-1\} $ yields
\begin{align*}
    \sum_{t=s_k}^{s_{k+1}-1}\! \left\{c_t(x_t, u_t) -c_t^{(I_t)}\right\}&\leq (1-\rho_{\delta})\alpha \sum_{t=s_k}^{s_{k+1}-1} \rho_{\delta}^{t-s_k}\leq \alpha.
\end{align*}
We then take the sum with respect to $ k $.
The initial interval contributes zero because both trajectories start from $x_0$ and use the same controller. Hence, the first term in~\eqref{decomposition} is at most $\alpha\bbE[N_I]$.
Moreover, by $N_I\le N_J$ and Proposition~\ref{theorem:switch}, this is upper bounded by $\alpha \eta \Delta_{\mathrm{c}} T / 4$.

We next bound the second term in \eqref{decomposition}.
Let $S_t^{(J)}:=\mathbb{I}[J_t\neq J_{t-1}]$ for $t\geq1$
indicate whether the tentative selection switches at time $t$.
By Lemma~\ref{lem:vbo}, we have
$
c_t^{(I_t)}-c_t^{(J_t)}
\leq \Delta_{\mathrm{c}}\mathbb{I}[I_t\neq J_t]
$.
Fix $t\geq1$ such that $I_t\neq J_t$, and let
$
s:=\max\left(
\{0\}\cup
\{r\in\{1,\ldots,t\}:I_r\neq I_{r-1}\}
\right)
$
be the most recent switching time of the applied controller,
where $s=0$ if no switch has occurred.

Since $I_t\neq J_t$, the dwell-time filter must be
in its locked period, so that $t-s<\tau$.
Moreover, $I_s=J_s$ holds by the filtering rule
if $s\geq1$, and by initialization if $s=0$.
Since the applied controller does not switch after $s$,
we also have $I_t=I_s=J_s$.
Thus, $I_t\neq J_t$ implies $J_s\neq J_t$,
which requires at least one switch in the tentative
selection during $\{s+1,\ldots,t\}$.
Consequently,
$
\mathbb{I}[I_t\neq J_t]
\leq \sum_{r=s+1}^{t}S_r^{(J)}\leq \sum_{r=\max\{1,t-\tau+1\}}^{t}S_r^{(J)}$.
The final bound on $\mathbb{I}[I_t\neq J_t]$ also holds when $I_t=J_t$ because its left-hand side is zero.
Thus, it follows from $
c_t^{(I_t)}-c_t^{(J_t)}
\leq \Delta_{\mathrm{c}}\mathbb{I}[I_t\neq J_t]
$, $I_0=J_0$, and Proposition~\ref{theorem:switch} that
\begin{equation*}
\begin{aligned}
\mathbb{E}\! \left[
\sum_{t=0}^{T-1}
\left(c_t^{(I_t)}-c_t^{(J_t)}\right)
\right]\leq \Delta_{\mathrm{c}}\sum_{t=1}^{T-1}
\sum_{r=\max\{1,t-\tau+1\}}^t
\mathbb{E}[S_r^{(J)}] \leq
\frac{\eta \Delta_{\mathrm{c}}^2}{4}
\sum_{t=1}^{T-1}\min\{t,\tau\} \leq \frac{\eta\tau \Delta_{\mathrm{c}}^2 T}{4} .
\end{aligned}
\end{equation*}

Lastly, by Proposition~\ref{theorem:switch}, the remaining term in \eqref{decomposition} can be bounded by the regret bound of Hedge (Proposition~\ref{pro:hedge}).
Combining the above bounds, we obtain \eqref{regretstandard}. Furthermore, letting $\eta=\Theta(\sqrt{T^{-1}\log N})$, we obtain $
        \regret_T = O\!\left(\sqrt{T \log N}\right)
    $. 
This completes the proof.
\end{proof}

The regret decomposition~\eqref{decomposition} shows that the control regret is bounded by the Hedge regret plus penalties for controller switching and dwell-time
filtering. This reveals a connection to the experts problem with switching costs~\citep{geulenregret}.
\red{Although lazy sampling and the dwell-time filter suppress switching, Hedge still updates its probabilities at every time step. Lazy sampling preserves the Hedge marginals while limiting the expected number of tentative switches to \(O(\eta T)\). Incremental stability bounds the cumulative cost discrepancy due to each actual switch by a constant. Moreover, the counterfactual costs under the tentative and filtered selections can differ only within \(\tau\) steps following a tentative switch. Consequently, the additional regret due to switching and the dwell-time filter is \(O(\eta T)\) for fixed \(\tau\), preserving the \(O(\sqrt{T\log N})\) regret order of the standard experts problem.}
See Appendix~\ref{app:discussion} for
further discussion.

\section{Numerical Experiments}\label{sec:experiment}
We consider $x_{t+1}=Ax_t+Bu_t+w_t$ with $x_0=0$ and four
candidate controllers $u_t=-K_ix_t$, where
$
A=\begin{bmatrix}1.05&0.5\\0&1.05\end{bmatrix}$, 
$B=\begin{bmatrix}0\\1\end{bmatrix}$,
and $K_1 = [3.20, 2.58]$, $K_2 = [0.69,0.83]$, $K_3 = [3.50,2.70]$, $K_4 = [0.50,0.75]$.
All candidates are Schur stabilizing, but
$\rho((A-BK_2)(A-BK_1))>1$, so that arbitrary switching need not be stable. 
They satisfy E-ISS and $\delta$-ES with $(\kappa,\rho) = (\kappa_\delta,\rho_\delta)=(5.2,0.78)$ and $\beta = \kappa / (1-\rho) $.
We use $\tau=7$, which satisfies
$\kappa\rho^\tau<1$,
and independent disturbances uniformly distributed over
$\{w\in\mathbb R^2:\|w\|\le0.2\}$.
The disturbance and cost-parameter sequences are generated before the online experiments and remain fixed
throughout each run, independently of the controller selections.
All methods use the same sequences, and each horizon $T$ uses their first $T$ elements.

\begin{figure}[!t]
  \centering
  \begin{minipage}[t]{0.49\columnwidth}
    \vspace{0pt}
    \centering
    \includegraphics[width=0.9\linewidth]{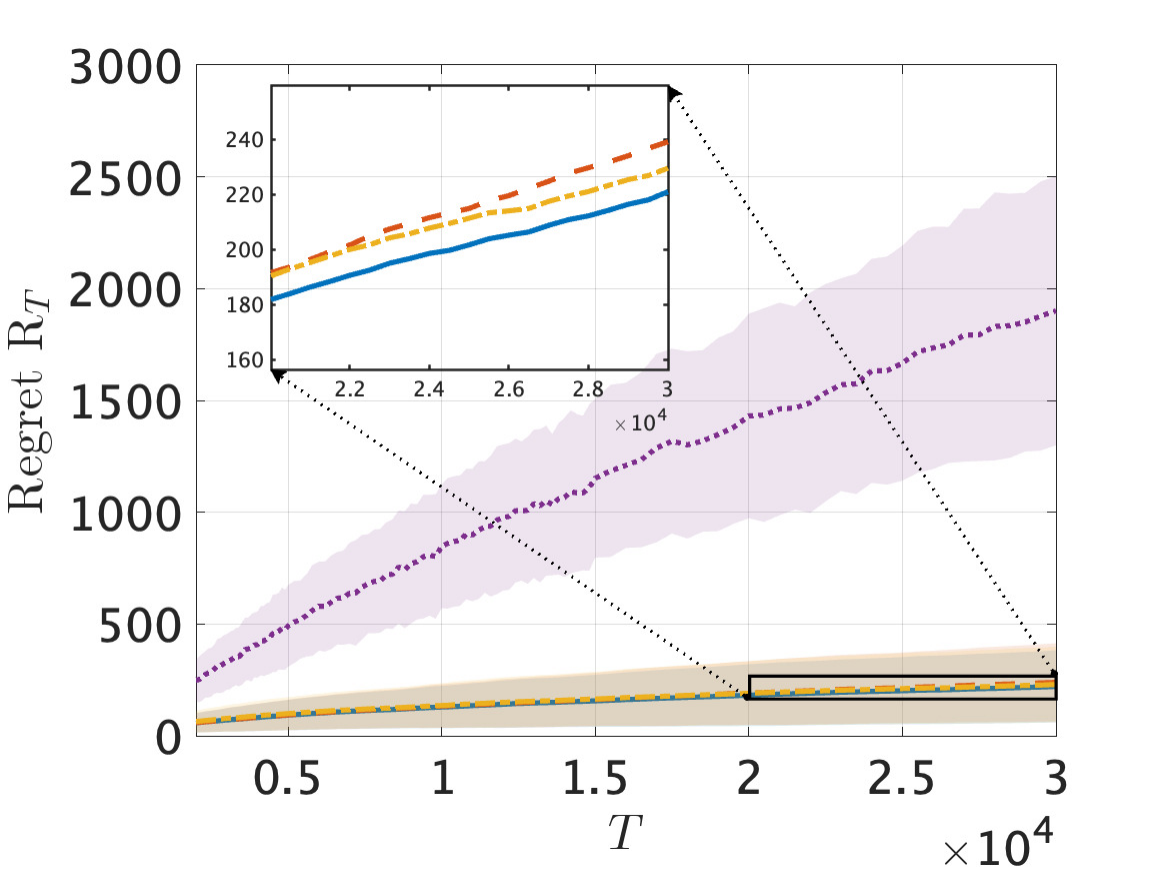}
    
    (a)
  \end{minipage}
  \hfill
  \begin{minipage}[t]{0.49\columnwidth}
    \vspace{0pt}
    \centering
    \includegraphics[width=0.82\linewidth]{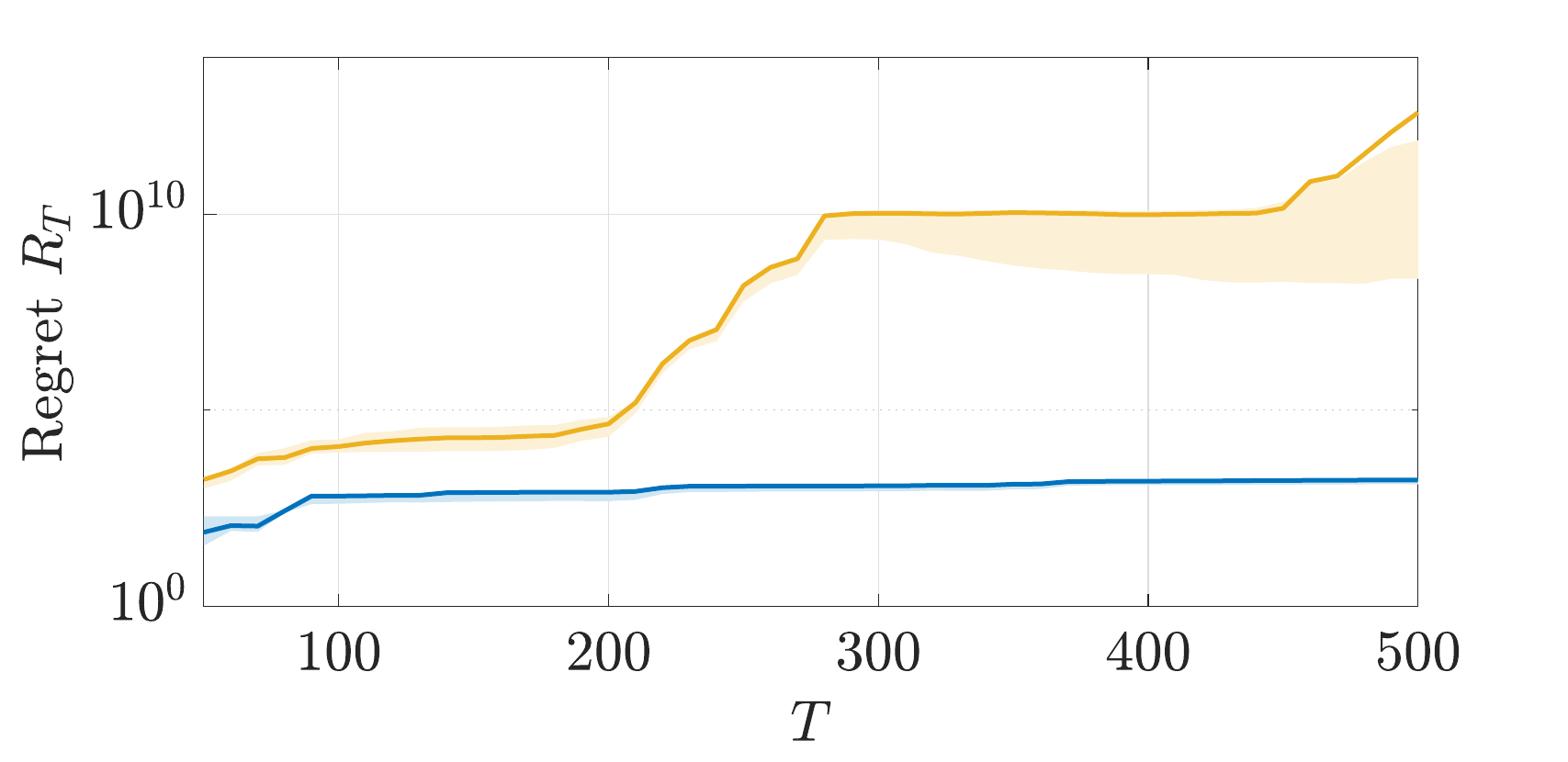}
    
    (b)
    
    \vspace{2mm}
    \includegraphics[width=0.82\linewidth]{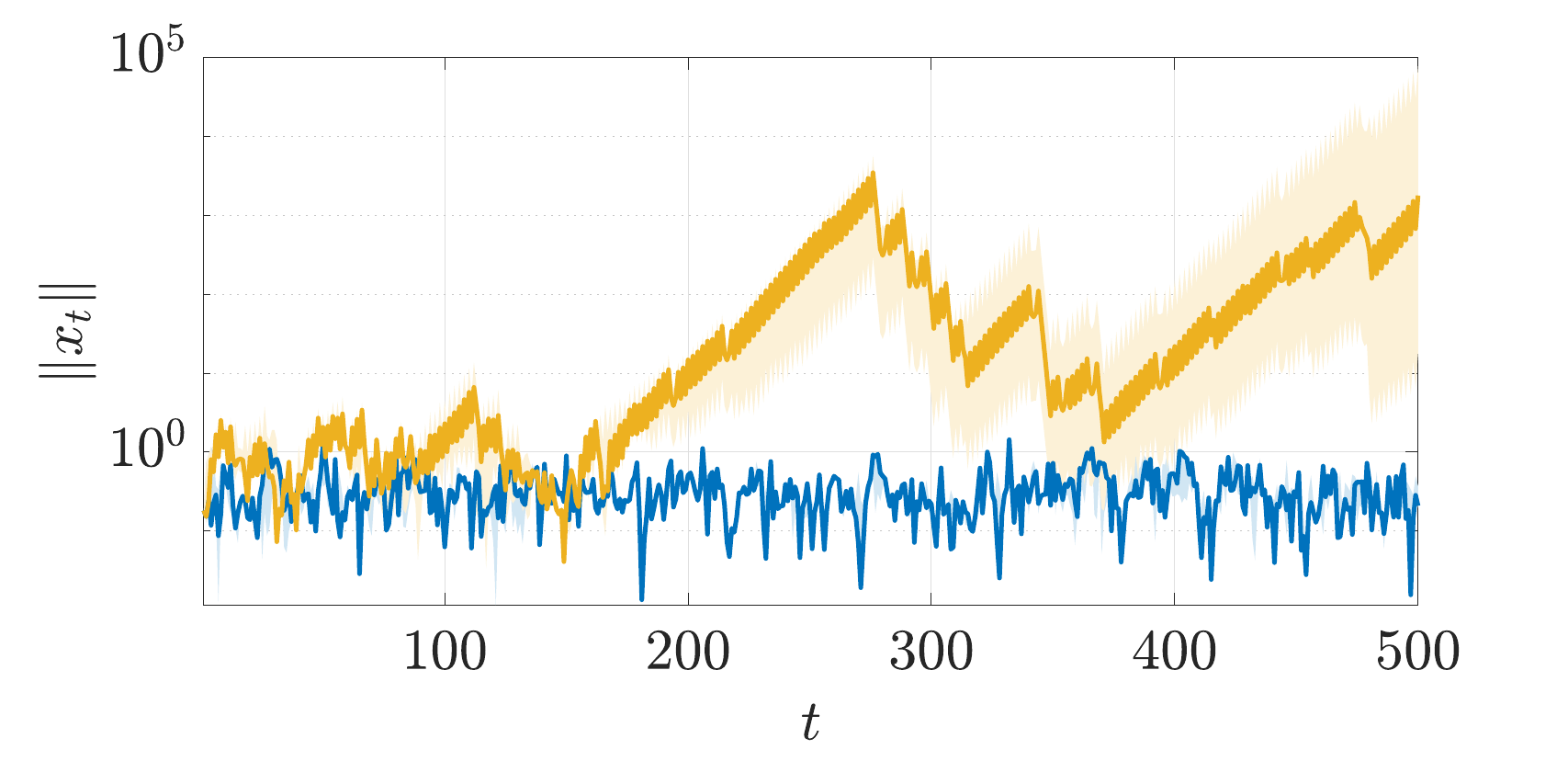}
    
    (c)
  \end{minipage}
  \caption{
  (a) Mean regret of the proposed method (solid),
  without lazy sampling (dashed),
  without the dwell-time filter (dash-dotted),
  and BAPS (dotted).
  Shaded regions indicate $\pm1$ standard deviation
  over $500$ trials.
  (b)--(c) Comparison of the proposed method (blue)
  and its variant without the dwell-time filter (yellow)
  over $300$ trials.
  (b) Mean regret.
  (c) State norm for $T=500$; lines \red{show the median, and shaded regions indicate} the
  $10$--$90$ percentile range.
  }
  \label{fig:regret_statenorm}
\end{figure}

We first compare the proposed method with independent Hedge sampling
(without lazy sampling), lazy sampling without the dwell-time filter,
and a bandit-feedback baseline, Bandit-based Adaptive Policy Selection (BAPS) \citep{anand2022online}. 
BAPS achieves a regret bound of $\widetilde{O}(T^{2/3}N^{1/3})$ in the bandit setting.  
We use
$c_t(x,u)=x^\top Q_tx+0.01u^2$, where $Q_t$ switches between
$\operatorname{diag}(43,1)$ and $\operatorname{diag}(1,3)$
with probability $0.25$ per step (initial probabilities $1/2$).
Following Theorem~\ref{theorem: regret}, we use $\eta=100\sqrt{T^{-1}\log(4)}$ for the proposed
method and the two variants without lazy sampling
and without the dwell-time filter.
For BAPS, we tune its hyperparameters by minimizing
the average normalized regret in independent pilot experiments (see Appendix~\ref{app:numerical_experiments} for details).
Fig.~\ref{fig:regret_statenorm}(a) shows that the proposed method achieves \red{much} lower mean regret than BAPS.
This verifies that learning performance is improved by utilizing the full-information feedback.
Although our regret guarantee does not apply to the two variants, they perform comparably in this example.

To examine the dwell-time filter, we construct
a time-varying cost sequence that encourages frequent
switching between controllers $K_1$ and $K_2$ (see Appendix~\ref{app:numerical_experiments} for details).
As shown in Figs.~\ref{fig:regret_statenorm}(b) and (c), removing the dwell-time filter results in \red{significantly} larger regret and state norms, whereas the proposed method achieves small regret while keeping the state norm small.
\red{This is because frequent switching between \(K_1\) and \(K_2\) can destabilize the system, whereas enforcing the prescribed minimum dwell time keeps the state uniformly bounded.}

\section{Conclusion}\label{sec:conclusion}
We proposed an online switching control algorithm for known nonlinear systems under full-information feedback. 
Lazy sampling preserves Hedge's marginals while controlling switching costs, and the dwell-time filter guarantees state boundedness. This two-stage mechanism yields an expected regret upper bound of $O(\sqrt{T\log N})$. Numerical examples illustrate the learning performance and show that removing the filter can amplify both state norms and regret. Future work includes extending this framework to unknown systems. In such scenarios, constructing counterfactual trajectories is no longer straightforward, and incorporating system identification may be needed.

\acks{This work was supported in part by JSPS KAKENHI Grant Numbers
JP24K17297, JP22H01508, and JP23K22778, and in part by JST under ASPIRE Grant Number JPMJAP2402.}

\bibliography{refs}

\clearpage
\appendix

\section{Proof of Lemma~\ref{lem:vbo}}\label{app:vbo}

Since $\{x_t^{(i)}\}$ is the state trajectory under the fixed
controller $\pi_i$, the E-ISS property in Assumption~\ref{ass:controller} yields
\[
\|x_t^{(i)}\|
\leq \kappa\rho^t\|x_0\|+\beta W
\leq \kappa\|x_0\|+\beta W
=X_c.
\]
Moreover, by the Lipschitz continuity of $\pi_i$,
\[
\begin{aligned}
\|u_t^{(i)}\|=\|\pi_i(x_t^{(i)})\|&\leq \|\pi_i(x_t^{(i)})-\pi_i(0)\|
     +\|\pi_i(0)\|\\
&\leq L_\pi\|x_t^{(i)}\|+\bar{\pi}_0\\
&\leq L_\pi X_c+\bar{\pi}_0=U_c.
\end{aligned}
\]

Next, for any $i,j\in\mathcal{P}$, the local Lipschitz
continuity of $c_t$ in Assumption~\ref{ass:cost} implies
\[
\begin{aligned}
|c_t^{(i)}-c_t^{(j)}|&\leq
\Bigl[
L_{c1}\bigl\{
\max(\|x_t^{(i)}\|,\|x_t^{(j)}\|)
+\max(\|u_t^{(i)}\|,\|u_t^{(j)}\|)
\bigr\}
+L_{c2}
\Bigr]\\
&\quad\times
\bigl(\|x_t^{(i)}-x_t^{(j)}\|
+\|u_t^{(i)}-u_t^{(j)}\|\bigr)\\
&\leq
\bigl\{L_{c1}(X_c+U_c)+L_{c2}\bigr\}
(2X_c+2U_c)\\
&=\Delta_{\mathrm c}.
\end{aligned}
\]
This completes the proof.
\hfill $\blacksquare$

\section{Proof of Proposition~\ref{theorem:switch}}\label{app:theorem:switch}

We first show that lazy sampling preserves the marginal
distributions, i.e., $\bbP[J_t=i]=p_{t,i}$ for all
$i\in\calP$ and $t\ge0$.
We prove the claim by induction on $t$.
The case $t=0$ holds by initialization.
Suppose that $\bbP[J_{t-1}=i]=p_{t-1,i}$
for all $i\in\calP$.
By Algorithm~\ref{alg:Lazy sampling}, the conditional probabilities satisfy
\[
\bbP[J_t=i\mid J_{t-1}=i]
=\min\! \left\{1,\frac{p_{t,i}}{p_{t-1,i}}\right\}
\]
and, for $j\neq i$,
\[
\begin{aligned}
&\bbP[J_t=i\mid J_{t-1}=j]=
\max\! \left\{0,1-\frac{p_{t,j}}{p_{t-1,j}}\right\}
\frac{(p_{t,i}-p_{t-1,i})_+}
{\sum_{k=1}^N(p_{t,k}-p_{t-1,k})_+},
\end{aligned}
\]
where $(a)_+:=\max\{0,a\}$.
Here, we first consider $p_t\neq p_{t-1}$,
so that the denominator is positive.

By the law of total probability,
\[
\begin{aligned}
\bbP[J_t=i]
&=\sum_{j=1}^N
\bbP[J_t=i\mid J_{t-1}=j]p_{t-1,j}\\
&=\min\{p_{t-1,i},p_{t,i}\}+
\frac{\sum_{j\neq i}(p_{t-1,j}-p_{t,j})_+}
{\sum_{k=1}^N(p_{t,k}-p_{t-1,k})_+}
(p_{t,i}-p_{t-1,i})_+.
\end{aligned}
\]
Since $p_t$ and $p_{t-1}$ are probability
distributions, we have
\[
\sum_{j=1}^N(p_{t-1,j}-p_{t,j})_+
=
\sum_{j=1}^N(p_{t,j}-p_{t-1,j})_+.
\]
If $p_{t,i}\ge p_{t-1,i}$, the second term
in the expression for $\bbP[J_t=i]$ equals
$p_{t,i}-p_{t-1,i}$.
If $p_{t,i}<p_{t-1,i}$, the second term
is zero.
Therefore, in either case,
\[
\begin{aligned}
\bbP[J_t=i]
&=\min\{p_{t-1,i},p_{t,i}\}+(p_{t,i}-p_{t-1,i})_+=p_{t,i}.
\end{aligned}
\]
If $p_t=p_{t-1}$, Algorithm~\ref{alg:Lazy sampling} maintains
the previous selection with probability one,
and hence the same conclusion holds.
Thus, the claim follows by induction.

The upper bound on the switching probability is established by combining the following two lemmas.
\ito{The first lemma characterizes the switching probability in terms of total variation distance.}
\begin{lemma}\label{lem:lazy-switch}
The sequence $ \{J_t\} $ generated by lazy sampling and marginal distributions $ \{p_t\} $
satisfies
\begin{equation*}
    \mathbb{P}[J_{t+1} \neq J_{t}] = \|p_{t+1} - p_t\|_{\mathrm{TV}},
\end{equation*}
where the total variation distance between $p_t$ and $p_{t+1}$ is defined by $\|p_{t+1}-p_t\|_{\mathrm{TV}} := \frac{1}{2} \sum_{i=1}^N |p_{t+1,i}-p_{t,i}|$.
\end{lemma}
\begin{proof}
The probability that the consecutive choices $ J_t $ and $ J_{t+1} $ coincide is given by
\begin{align*}
    \mathbb{P}[(J_{t+1}, J_t)=(i,i)]&=\mathbb{P}[J_{t+1}=i\mid J_t=i]\mathbb{P}[J_t=i]\\
    &=\min\!\left\{1,\frac{p_{t+1,i}}{p_{t,i}}\right\}p_{t,i}\\
    &=\min\{p_{t,i}, p_{t+1,i}\}. 
\end{align*}
Thus, the probability of retaining the same selection between times $t$ and $t+1$ is
\begin{align*}
    \mathbb{P}[J_{t+1}=J_t]
    &= \sum_{i=1}^N \min\{p_{t,i}, p_{t+1,i}\}\\
    &=\sum_{i=1}^N \frac{p_{t,i}+p_{t+1,i}-|p_{t+1, i}-p_{t,i}|}{2}\\
    &= 1-\frac{1}{2}\sum_{i=1}^N |p_{t+1, i}-p_{t,i}|,
\end{align*}
where we used $\min\{a,b\}=\frac{a+b-|b-a|}{2}$. 
The lemma immediately follows from $\mathbb{P}[J_{t+1}\neq J_t]=1-\mathbb{P}[J_{t+1}= J_t]$. 
\end{proof}

By the coupling inequality
\citep[Proposition~4.7]{levin2017markov}, if $ J_t \sim p_t $ and $ J_{t+1} \sim p_{t+1} $, then $ \bbP[J_{t+1} \neq J_t] \ge \| p_{t+1} - p_t \|_{\mathrm{TV}} $. Therefore, Lemma~\ref{lem:lazy-switch} shows that lazy sampling attains the maximum probability that $ J_t $ and $ J_{t+1} $ coincide subject to $ J_t \sim p_t $ and $ J_{t+1} \sim p_{t+1} $.

The next lemma bounds the total variation distance of distributions given by Hedge.
\begin{lemma}\label{lem:hedge-tv}
Let $p_{t,i}$ be defined by \eqref{hedge:update}
with a constant learning rate $\eta>0$, and let
$m_t:=\min_{i\in\calP}c_t^{(i)}$ and
$M_t:=\max_{i\in\calP}c_t^{(i)}$.
Then, it holds that
\begin{equation}
    \|p_{t+1} - p_{t}\|_{\mathrm{TV}}\leq \frac{\eta (M_t-m_t)}{4}.
\end{equation}
\end{lemma}
\begin{proof}
Since the learning rate is constant, \eqref{hedge:update} is equivalent to
the recursive form
\begin{equation}\label{eq:hedge_rec}
    p_{t+1,i}
    = \frac{p_{t,i} e^{-\eta c_{t}^{(i)}}}
           {\sum_{j \in \calP} p_{t,j} e^{-\eta c_{t}^{(j)}}} .
\end{equation}
Hence, the Kullback--Leibler divergence between $ p_t $ and $ p_{t+1} $ is bounded as
\begin{align}
    D_{\mathrm{KL}}(p_t \| p_{t+1})
    &:= \sum_{i \in \calP} p_{t,i} \log \frac{p_{t,i}}{p_{t+1,i}} \notag \\
    &= \sum_{i \in \calP} p_{t,i} \! \left\{ \eta c_{t}^{(i)}
    + \log \! \left( \sum_{j \in \calP} p_{t,j} e^{-\eta  c_{t}^{(j)} } \right)
    \right\} \notag \\
    &= \eta \sum_{i \in \calP} p_{t,i} c_{t}^{(i)}
    + \log \! \left( \sum_{j \in \calP} p_{t,j} e^{-\eta c_{t}^{(j)}} \right)
    \notag \\
    &\leq \frac{\eta^2 (M_t - m_t)^2}{8}, \label{eq:kl}
\end{align}
where the last inequality follows from Hoeffding's
lemma.
Pinsker's inequality then yields
\begin{equation*}
    \| p_{t+1} - p_t \|_{\mathrm{TV}}
    \leq \sqrt{\tfrac{1}{2} D_{\mathrm{KL}}(p_t \| p_{t+1})}
    \leq \frac{\eta (M_t - m_t)}{4} .
\end{equation*}
This completes the proof.
\end{proof}

Proposition~\ref{theorem:switch} follows by applying
Lemmas~\ref{lem:lazy-switch} and \ref{lem:hedge-tv}, together with $ M_t - m_t \leq \Delta_{\mathrm{c}} $ by Lemma~\ref{lem:vbo}.

\section{Proof of Proposition~\ref{th:statebound}}\label{app:statebound}

\ito{Let $0=t_0<t_1<\cdots<t_{N_I} <  T$ denote the actual switching times, where $t_0$ represents initialization, and set $t_{N_I + 1} := T$.
For $k\in \{0,\ldots,N_I-1\}$, the dwell-time constraint
ensures $t_{k+1}-t_k\ge\tau$. Thus, E-ISS inequality \eqref{eq:E-ISS} yields that
\begin{equation}\label{eq:bound_tk1}
    \| x_{t_{k+1}} \| \le \kappa \rho^{t_{k+1} - t_k} \| x_{t_k} \| + \beta W \le \kappa \rho^\tau \| x_{t_k} \| + \beta W , \qquad \forall k\in \{0,\ldots,N_I-1\} .
\end{equation}
Moreover, by the assumption $ \tau > \frac{\log \kappa}{\log (1/\rho)} $, we have $ \gamma := \kappa\rho^\tau < 1 $.
Thus, it follows from \eqref{eq:bound_tk1} that
\begin{equation*}
    \| x_{t_k} \| \le \gamma^k \| x_0 \| + \sum_{s=1}^k \gamma^{s-1} \beta W \le \| x_0 \| + \frac{\beta W}{1-\gamma} , \qquad \forall k\in \{0,\ldots,N_I\} .
\end{equation*}
For each $t\in\{0,\ldots,T\}$, let
$k=\max\{j\in\{0,\ldots,N_I\}:t_j\le t\}$.
Applying E-ISS from $t_k$ yields
\begin{align*}
    \| x_t \| \le \kappa \rho^{t-t_k} \| x_{t_k} \| + \beta W  &\le \kappa \| x_{t_k} \| + \beta W \le \kappa \left( \| x_0 \| + \frac{\beta W}{1-\gamma}  \right) + \beta W ,
\end{align*}
which completes the proof.
\hfill $\blacksquare$}

\section{Discussion}\label{app:discussion}
The proof of Theorem~\ref{theorem: regret} implies that the online switching control problem addressed in this paper reduces to the experts problem with switching costs. 
The impact of switching costs differs significantly between full-information and bandit settings. 
In the full-information setting, the optimal regret remains $\Theta(\sqrt{T\log N})$ regardless of the presence of switching costs~\citep{littlestone1994weighted, geulenregret}.
However, in the bandit setting, the presence of switching costs degrades the optimal regret from $\Theta(\sqrt{TN})$ to $\wtilde{\Theta}(T^{2/3}N^{1/3})$~\citep{dekel2014bandits}.

Previous bandit-feedback approaches use batching to ensure
a minimum dwell time and reduce transient
effects in performance evaluation~\citep{anand2022online,li2023online}.
In our full-information setting, counterfactual costs
eliminate the latter need.
Consider batched Hedge with batch length $b\ge\tau$
and independent sampling at batch boundaries.
Then, it can be shown that its regret satisfies
\begin{equation*}
    \regret_T \leq \frac{\alpha T}{b}
    + \frac{\log N}{\eta}
    + \frac{\eta\Delta_{\mathrm c}^2 bT}{8}.
\end{equation*}
Optimizing this bound over $b$ and $\eta$ for independently sampled batched Hedge yields
$O(T^{2/3}(\log N)^{1/3})$, whereas our approach
achieves $O(\sqrt{T\log N})$ through lazy sampling.

\section{Additional Details on the Numerical Experiments}
\label{app:numerical_experiments}

\subsection{Regret Comparison and Parameter Selection}

For BAPS, we use Algorithm~2 of
\citet{anand2022online}.
The learning-rate and batch-length coefficients
were tuned by grid search over
\[
\begin{aligned}
\mathcal{G}_{\eta}
&=\{0.001,0.01,0.03,0.1,0.3,1,3,10,30,100\},\\
\mathcal{G}_{b}
&=\{0.25,0.5,1,2,4,8,16,32\}.
\end{aligned}
\]
We evaluated each of the $80$ parameter combinations
using $24$ independent trials at each horizon
$T\in\{3000,12000,30000\}$.
The parameters were selected to minimize the average
terminal regret normalized by $T^{2/3}$ across these
horizons.
The resulting learning rate and batch length are
\[
\eta_{\mathrm B}
=\frac{3}{T^{2/3}4^{1/3}},
\qquad
b_{\mathrm B}
=\max\{7,\lceil 2(T/4)^{1/3}\rceil\}.
\]

\subsection{Adversarial Cost Construction for the Dwell-time Stress Test}
In the second experiment, we construct a quadratic cost sequence,
\begin{equation*}
c_t(x,u)=x^\top Q_tx+0.01u^2,\qquad
Q_t=0.1I+\lambda_t v_tv_t^\top,
\end{equation*}
where $\|v_t\|_2=1$ and $0\le\lambda_t\le2000$.
The construction is intended to encourage frequent switching between
$K_1$ and $K_2$, without prescribing the controller selections.
We first sample one disturbance sequence $w_0,\ldots,w_{499}$ and
compute the four fixed-controller counterfactual trajectories
$(x_t^{(i)},u_t^{(i)})$. Initially, let $Q_0=Q_1=0.1I$.
For $t\ge2$, define
\begin{equation*}
\begin{aligned}
S_t&:=x_t^{(1)}x_t^{(1)\top}-x_t^{(2)}x_t^{(2)\top},\\
g_t^0&:=0.1\operatorname{tr}(S_t)
 +0.01\bigl((u_t^{(1)})^2-(u_t^{(2)})^2\bigr),\\
D_{t-1}&:=\sum_{s=0}^{t-1}
 \bigl(c_s(x_s^{(1)},u_s^{(1)})
       -c_s(x_s^{(2)},u_s^{(2)})\bigr),\\
\zeta_t&:=(-1)^{t+1}-D_{t-1}-g_t^0.
\end{aligned}
\end{equation*}
To explain this construction, let
$d_t:=c_t^{(1)}-c_t^{(2)}$.
By the definition of $Q_t$, we have
$
d_t=g_t^0+\lambda_t v_t^\top S_t v_t
$,
and hence $D_t=D_{t-1}+d_t$.
Since the ratio of their Hedge selection probabilities is
$p_{t+1,1}/p_{t+1,2}=\exp(-\eta D_t)$,
we aim to make $D_t$ alternate between $-1$ and $1$ to encourage frequent changes in their relative
selection probabilities.
The quantity $\zeta_t$ represents the additional cost
difference required to achieve the target
$D_t=(-1)^{t+1}$, namely,
$
\lambda_t v_t^\top S_t v_t=\zeta_t
$.
We therefore choose $v_t$ along an extremal eigenvector
of $S_t$ to adjust the cost difference in the desired
direction, and determine $\lambda_t$ subject to its
upper bound.

If $\zeta_t>0$, choose $\mu_t$ as the largest eigenvalue of
$S_t$; if $\zeta_t<0$, choose its smallest eigenvalue.
Let $v_t$ be a corresponding unit eigenvector. Since $v_t^\top S_tv_t=\mu_t$, the target
$\lambda_t v_t^\top S_tv_t=\zeta_t$ is achieved by
$\lambda_t=\zeta_t/\mu_t$, provided that
$\zeta_t\mu_t>0$.
To satisfy $0\le\lambda_t\le2000$ and avoid division
by a nearly zero eigenvalue, we set
\begin{equation*}
\lambda_t=
\begin{cases}
\min\{2000,|\zeta_t/\mu_t|\},
& \zeta_t\mu_t>0,\ |\mu_t|>10^{-12},\\
0, & \text{otherwise}.
\end{cases}
\end{equation*}
If $\zeta_t=0$, we also set $\lambda_t=0$.
Using the resulting $Q_t$, we compute the
counterfactual costs and update $D_t$.
Although the constraints on $\lambda_t$ may prevent
the alternating target gap from being attained exactly,
$0.1I\preceq Q_t\preceq2000.1I$ holds for all $t$.
Thus, Assumption~\ref{ass:cost} holds.
Since the cost sequence is constructed entirely from
presampled disturbances and fixed-controller
counterfactual trajectories, it is independent of
the online controller selections and thus satisfies
Assumption~\ref{adverse}.

\end{document}